\documentclass[12pt,a4paper]{article}
\usepackage{amssymb}
\newtheorem{dfn}{Definition}[section]
\newtheorem{tetel}[dfn]{Theorem}

\newtheorem{lemma}[dfn]{Lemma}

\newcommand{\R}{\mathbb R}
\newcommand{\C}{\mathbb C}
\newcommand{\N}{\mathbb N}
\newcommand{\Z}{\mathbb Z}
\newcommand{\T}{\mathbb T}
\def\Re{\hbox{Re}\,}
\def\Tr{\mbox{Tr}\,}
\def\Prob{\hbox{Prob}\,}
\def\Diag{{\rm Diag}}
\def\<{\langle}
\def\>{\rangle}
\def\Im{\hbox{Im}\,}
\newcommand{\ud}{\,d}

\newcommand{\tr}{\mathop{\mathrm{Tr}}}
\newcommand{\vv}{random variable}

\newcommand{\bv}{\nobreak\hfill $\square$}
\def\im{{\rm i}}
\begin{document}
\vskip 1cm
\centerline{\LARGE{\bf On asymptotics of large}}
\medskip
\centerline{\LARGE{\bf Haar distributed unitary matrices}\footnote{This work
was partially supported by OTKA T032662 and T032374.}}
\medskip \bigskip
\centerline{D\'enes Petz\footnote{E-mail: petz@math.bme.hu} 
and J\'ulia R\'effy\footnote{E-mail: reffyj@math.bme.hu}}
\bigskip
\medskip
\centerline{Department for Mathematical Analysis}
\centerline{Budapest University of Technology and Economics}
\centerline{ H-1521 Budapest XI., Hungary}
\bigskip\bigskip
\noindent
\begin{quote}
Let $U_n$ be an $n \times n$ Haar unitary matrix. In this paper,
the asymptotic normality and independence of $\Tr U_n, \Tr U_n^2,
\dots, \Tr U_n^k$ are shown by using elementary methods. More
generally, it is shown that the renormalized truncated Haar
unitaries converge to a Gaussian random matrix in distribution. 
\end{quote}

\begin{quote}
\end{quote}
\bigskip\bigskip
\noindent

\section{Introduction}
Entries of a random matrix are random variables but a random matrix is 
equivalently considered as a probability measure on the set of matrices.
A simple example of random matrix has independent identically distributed
entries. In this paper random unitary matrices are studied whose
entries must be correlated.
  
A  unitary matrix $U=(U_{ij})$ is a matrix with complex entries and
$UU^*=U^*U=I$. In terms the entries these relations mean that
\begin{eqnarray}
&& \sum_{j=1}^n |U_{ij}|^2=\sum_{i=1}^n |U_{ij}|^2=1, \hbox{for all\ } 1\le
i,j\le n, \\ &&
\sum_{l=1}^n U_{il}\overline U_{jl}, \hbox{for all\ } 1\le i,j \le n,
\quad i\ne j.
\end{eqnarray}

The set $\mathcal U(n)$ of $n\times n$ unitary matrices forms a compact 
topological group with respect to the matrix multiplication and the usual
topology, therefore there exists a unique (up to the scalar
multiplication) translation invariant measure on $\mathcal U(n)$, the
so-called Haar measure. We will consider a random variable $U_n$ 
which maps from a probability space to $\mathcal U(n)$, and take 
its values uniformly from $\mathcal U(n)$, i.e. if $H\subset \mathcal U(n)$,
then 
\[
\Prob(U_n\in H)=\gamma (H),
\]
where $\gamma$ is the normalized Haar measure on $\mathcal U(n)$.
We call this random variable a Haar unitary random variable, or shortly Haar unitary. 

Relations (1) and (2) show that the column vectors of a unitary are
pairwise othogonal unit vectors and the distribution of each column 
has unitarily invariant joint distribution in case of a Haar distributed
unitary. This fact allows a simple construction of a Haar unitary from
a Gaussian matrix with i.i.d. entries. In Sect. 2 this construction is
written up in details. The construction allows to determine the 
distribution of the matrix elements. Due to the permutation invariance,
the elements are identically distributed. We observe that large powers
of the eigenvalues are independent and uniformly distributed. The main
aim of the paper is to
study asymptotic questions when the matrix size is going to infinity.
In this limit the matrix elements are going to be Gaussian (after a 
renormalization) and more generally, the truncated matrix converges to
a Gaussian matrix. This is the content of Sect. 3. In the rest of the 
paper we show that the trace of the powers is going to Gaussian in
the limit, moreover the traces of different powers are asymptotically
independent.  Actually, this has been shown by Diaconis and Shahshahani
\cite{diaconis}, they determined the Fourier transform of the limit 
distribution of the eigenvalues. In their proof the characters of
the symmetric and unitary groups, and different bases of the symmetric
polynomials play the main role. Here we get the same result in a more 
elementary way. The method of moments is used, and we examine the
order of magnitude of the different summands in the trace. This method
could be familiar from Arnold \cite{arnold} when he studied Wigner
matrices, or from Bai and Yin \cite{bai} for sample covariance
matrices. 

Random unitary matrices may be applied to several physical phenomena,
such as chaotic scattering or statistical properties of periodically
driven quantum systems \cite{Haake}. The large deviation theorem for
the emprical eigenvalue density of Haar unitaries (and of some other
unitaries) was established in \cite{HPuni}. Independent Haar unitaries
provide also a simple example of asymptotical freeness \cite{HP, Voi}.
 
\section{Haar unitary matrices}\label{elso}

Let $\xi$ be a complex-valued random variable. If $\Re \xi$ and $\Im
\xi$ are independent and normally distributed according to $N(0,1/2)$, 
then we call $\xi$ a {\it standard complex normal variable}. The
terminology is justified by the properties $E(\xi)=0$ and $E(\xi
\overline{\xi})=1$.

\begin{lemma}\label{normom}
Assume that $R \ge 0$ and $R^2$ has exponential distribution with parameter 1,
$\theta$ is uniform on the interval $[0,2\pi]$, and assume that $R$ and
$\theta$ are independent. Then $\xi =R e^{i\theta}$ is a standard
complex normal {\vv}
and
\[
E(\xi^k\overline \xi^\ell)=\delta_{k\ell}k! \qquad (k,\ell \in \Z_+)
\]
\end{lemma}

\begin{biz}
Let $X$ and $Y$ be real-valued random variables and assume that $X+iY$ is 
standard complex normal. For $r>0$ and $0\le \theta_0 \le 2\pi$
set 
\[S_{r,\theta_0}:=\{\rho e^{\im\psi}:0\le \rho\le r,0\le \psi\le\theta_0\},\] then
\begin{eqnarray*}
P\left(X+\im Y\in S_{r,\theta_0}\right)&&\kern-.24in
=\frac{1}{\pi}\int\int_{\{(s,t):s+\im t\in S_{r,\theta_0}\}}e^{-\left(s^2+t^2\right)}
dsdt\\
&&\kern-.24in =\frac{1}{\pi}\int_0^{\theta_0}d\psi\int_0^r\rho e^{-\rho^2}d\rho\\
&&\kern-.24in = \frac{1}{2\pi}\theta_0\left(1-e^{-r^2}\right)=P\left(\xi \in S_{r,\theta_0}\right).
\end{eqnarray*}
This
proves the first part which makes easy to compute the moments:
\[
E(\xi^k\overline \xi^\ell)=E\left(R^{k+\ell}\right)E\left(e^{i\theta(k-\ell)}\right)=
\delta_{k\ell}E(R^{2k}),
\]
so we need the moments of the exponentional distribution. We have by
partial integration
\begin{eqnarray*}
\int_0^\infty x^k e^{- x} \ud x& = &
-\left[x^ke^{- x}\right]_0^{\infty}+k\int_0^\infty x^{k-1}e^{-x}\ud x=
k\int_0^\infty x^{k-1}e^{-x}\ud x \\
&=&k(k-1)\int_0^\infty x^{k-2}e^{-x}\ud x=\dots=k!\int_0^\infty e^{-x}\ud x=k!
\end{eqnarray*}
which completes the proof of the lemma.
\end{biz}\bv

Next we recall how to get a Haar unitary from a Gaussian matrix with
independent entries by the Gram-Schmidt orthogonalization procedure on 
the column vectors.

Suppose that we have a complex random matrix $Z$ whose entries 
$Z_{ij}$ are mutually independent standard complex normal \vv s. 
We perform the Gram-Schmidt ortogonalization procedure on the column
vectors ${Z}_i$ ($i=1,2,\dots,n$), i.e. 
\begin{equation}
{U}_i=\frac{{Z}_i-\displaystyle{\sum_{l=1}^{i-1}} \< Z_i,U_l\>U_l}
{\left\| {Z}_i-\displaystyle{\sum_{l=1}^{i-1}} \< Z_i,U_l\>{U}_l \right\| }
,
\end{equation}
where
$$
\| (X_1,X_2,\dots,X_n)\|=\sqrt{\sum_{k=1}^n |X_k|^2}\,.
$$

\begin{lemma}
The above column vectors ${U_i}$ constitute a unitary matrix
${U}=\left({U_i}\right)_{i=1,\dots, n}$. Moreover, 
for all $V \in \mathcal U(n)$ the distributions of 
${U}$ and $V{U}$ are the same. 
\end{lemma}

\begin{biz}
First we prove, that for any $V \in \mathcal U(n)$ the matrices  ${Z}$
and $V{Z}$ have the same distribution. The entries 
$\xi_{ij}$ of $V Z $ are standard complex normal. Indeed,
$$
{\xi_{ij}=\sum_{l=1}^nV_{il}Z_{lj}}
$$
is normal. Furthermore
\[
E(\xi_{ij})=\sum_{l=1}^nV_{il}E(Z_{lj})=0,
\]
and
\[
E(\xi_{ij}\overline{\xi}_{ij})=\sum_{l=1}^n |V_{il}|^2E(Z_{lj}\overline 
Z_{lj})=\sum_{l=1}^n|V_{il}|^2=1\,.
\]

Next we prove that the correlation between two entries is zero. (In the
case of normally distributed \vv s this is equivalent to the independence.)
\begin{eqnarray*}
\lefteqn{E(\xi_{ij}\overline{\xi}_{sr})=
E\left(\left(\sum_{l=1}^n V_{il}Z_{lj}\right)\left(\sum_{k=1}^n
\overline{V_{sk}Z_{kr}}\right)\right)}\\
&&=\sum_{l=1}^n\sum_{k=1}^n
V_{il}\overline V_{sk}E(Z_{lj}\overline Z_{kr})
=\delta_{jr}\sum_{l=1}^nV_{il}\overline V_{sl}=\delta_{jr}\delta_{is}.
\end{eqnarray*}

The $i$th column of $VU$ is exactly $VU_i$ and we have
\begin{equation}
V{U}_i=\frac{V{Z}_i-\displaystyle{\sum_{l=1}^{i-1}} \< Z_i,U_l\>VU_l}
{\left\| {Z}_i-\displaystyle{\sum_{l=1}^{i-1}} \< Z_i,U_l\>{U}_l \right\| }
=
\frac{V{Z}_i-\displaystyle{\sum_{l=1}^{i-1}} \< VZ_i,VU_l\>VU_l}
{\left\|V {Z}_i-\displaystyle{\sum_{l=1}^{i-1}} \<V Z_i,VU_l\>V{U}_l \right\| }
\end{equation}
which is the Gram-Schmidt ortogonalization of the vectors $V{Z_i}$.
Since we showed above that $Z$ and $VZ$ are identically distributed, we 
conclude that $U$ and $VU$ are identically distributed as well.
Since the left invariance characterizes the Haar measure  on a compact
group, the above constructed ${U}$ is Haar distributed and its
distribution is right invariant as well.
\end{biz}\bv

The column vectors of a unitary matrix are pairwise orthogonal unit
vectors. On the bases of this fact we can determine a Haar unitary in
a slightly different way. The complex unit vectors form a compact
space on which the unitary group acts transitively. Therefore, there
exist a unique probability measure invariant under the action. Let us
call this measure uniform. To determine a Haar unitary, we choose
the first column vector $U_1$ uniformly from the space of $n$-vectors.
$U_2$ should be taken from the $n-1$ dimensional subspace orthogonal to
$U_1$ and choose it uniformly again. In general, if already $U_1,U_2,
\dots, U_j$ is chosen, we take $U_{j+1}$ from the $n-j$ dimensional
subspace orthogonal to $U_1,U_2, \dots, U_j$, again uniformly. The
column vectors constitute a unitary matrix and we check that its 
distribution is left invariant. Let $V$ be a fixed unitary. We show
that the vectors $VU_1, VU_2,\dots, VU_n$ are produced by the above
described procedure. They are obviously pairwise othogonal unit vectors.
$VU_1$ is uniformly distributed by the invariance property of the
distribution of $U_1$. Let $V(1)$ be such a unitary that
$V(1)VU_1=VU_1$. Then $V^{-1}V(1)VU_1=U_1$ and the choice of $U_2$
gives that $V^{-1}V(1)VU_2 \sim U_2$. It follows that
$V(1)VU_2 \sim VU_2$. Since $V(1)$ was arbitrary $VU_2$ is uniformly
distributed in the subspace orthogonal to $VU_1$. Similar argument works
for $VU_3, \dots, VU_n$. The Gram-Schmidt orthogonalization of the
columns of a Gaussian matrix gives a concrete realization of this 
procedure.

The permutation matrices are in $\mathcal U(n)$, and by multiplying 
with an appropriate permutation matrix every row and column can be 
transformed to any other row or column, so the translation invariance of 
a Haar unitary $U$ implies that all the entries have the same
distribution. From the above construction of a Haar unitary one can 
deduce easily the distribution of the entries:
$$
\frac{n-1}{\pi}(1-r^2)^{n-2}r \,dr\,d\theta
$$
(see also p. 140 in \cite{HP}). Since
$$
P(|\sqrt{n}U_{ij}|^2\ge x)=\Big(1-\frac{x}{n}\Big)^{n-1}\to e^{-x}
$$
$\sqrt{n}U_{ij}$ converges to a standard complex normal variable.
The correlation coefficient between $|U_{ii}|^2$ and $|U_{jj}|^2$ is
$1/(n-1)^2$ if $i\ne j$ (see p. 139 in \cite{HP}). 

In the next section we need the following technical lemma which tells us
that the expectation of many product of the entries are vanishing.

\begin{lemma}\label{mom0}{\bf (\cite{HP})}
Let  $i_1, \dots, i_h, j_1,\dots j_h \in
\{1,\dots ,n\}$ and $k_1,\dots ,k_h,m_1,\dots ,m_h$ be positive integers
for some $h \in \N$. If  
$$
{\sum_{i_r=u}(k_r-m_r)\ne 0} \qquad
\textrm{for some}\qquad 1\le u\le n$$
or
$$
{\sum_{j_r=v}(k_r-m_r)\ne 0}\qquad \textrm{for some}
\qquad 1\le v\le n,,
$$ 
then
\[
E\left((U_{i_1j_1}^{k_1}\overline U_{i_1j_1}^{m_1})\dots (U_{i_hj_h}^{k_h}
\overline U_{i_hj_h}^{m_h})\right)=0.
\]
\end{lemma}

\begin{biz}
Suppose that $t:=\sum_{i_r=u}(k_r-m_r)\ne 0$. 
The translation invariance of ${U}$ implies that multiplying this matrix by
$V={\Diag}(1,\dots, 1,e^{i\theta},1,\dots,1)\in \mathcal U(n)$ from the 
left  we get 
\[
E\left((U_{i_1j_1}^{k_1}\overline U_{i_1j_1}^{m_1})\dots (U_{i_hj_h}^{k_h}
\overline U_{i_hj_h}^{m_h})\right)=e^{it\theta}
E\left((U_{i_1j_1}^{k_1}\overline U_{i_1j_1}^{m_1})\dots (U_{i_hj_h}^{k_h}
\overline U_{i_hj_h}^{m_h})\right),
\]
for all $\theta \in \R$. 
\end{biz}\bv

Let $U$ be a Haar distributed $n \times n$ unitary matrix with
eigenvalues $\lambda_0,\lambda_1,\dots, \lambda_{n-1}$. The eigenvalues are
random variables with values in $\T:=\{z\in \C: |z|=1\}$, their joint 
distribution is well-known:
\begin{equation}
\frac{1}{n!}\prod_{i<j} |z_i-z_j|^2= \frac{1}{n!}\prod_{i<j} |e^{\im \theta_i}
-e^{\im \theta_j}|^2
\end{equation}
with respect to $dz_0\,dz_1\dots\, dz_{n-1}$, where $dz_i=d\theta_i/2\pi$
for $z_i=e^{\im \theta_i}$ (see p. 135 in \cite{HP}).

\begin{tetel}
For $m > n$ the random variables $\lambda_0^m,\lambda_1^m,\dots, 
\lambda_{n-1}^m$ are independent and uniformly distributed on $\T$.
\end{tetel}

\begin{biz}
Since the Fourier transform determines the joint distribution measure
of $\lambda_0^m,\lambda_1^m,\allowbreak \dots, \lambda_{n-1}^m$
 uniquely, it suffices to show that  
\begin{equation}\label{E:int}
\int z_0^{k_0 m}z_1^{k_1 m}\dots z_{n-1}^{k_{n-1}m} 
\prod_{i<j} |z_i-z_j|^2\,dz =0
\end{equation}
if at least one $k_j \in \Z$ is different from 0 ($dz=dz_0\,dz_1\dots 
dz_{n-1}$ and integration is over $\T^n$).
 
Let
\begin{equation}
\Delta(z_0 ,z_1, \dots, z_{n-1}):=\prod_{i<j}(z_i-z_j)=
\det\bigl[z^k_i\bigr]_{0\le i\le n-1,\,0\le k\le n-1}.
\end{equation}
(What we have here is the so-called Vandermonde determinant.) Then
\[
\prod_{i<j} |z_i-z_j|^2 = \Delta(z_0 ,z_1, \dots, z_{n-1})
\Delta(z_0^{-1} ,z_1^{-1}, \dots, z_{n-1}^{-1})\]
and one can write (\ref{E:int}) as an $n$-times complex contour integral
along the positively oriented $\T$:
\begin{eqnarray*}
\lefteqn{\int z_0^{k_0 m}z_1^{k_1 m}\dots 
z_{n-1}^{k_{n-1}m} 
\Delta(z_0 ,z_1, \dots, z_{n-1})\Delta(z_0^{-1} ,z_1^{-1}, \dots, z_{n-1}^{-1}) \,dz}\\
&&\kern-.2in= \int z_0^{k_0 m}\dots 
z_{n-1}^{k_{n-1}m} \sum_{\pi\in S_n}(-1)^{\sigma(\pi)}
z_0^{\pi(0)}\dots z_{n-1}^{\pi(n-1)} \sum_{\rho\in S_n}(-1)^{\sigma(\rho)}
z_0^{-\rho(0)}\dots z_{n-1}^{-\rho(n-1)}\,dz
\end{eqnarray*}
which is calculated according to the theorem of residue.
We need to find the coefficient of $z_0^{-1}z_1^{-1}\dots z_{n-1}^{-1}$, so we are looking for
the permutations, for which
$k_jm+\pi(j)-\rho(j)=-1$ if $0\le j \le n-1$, so $k_jm=\rho(j)-\pi(j)-1$. Here
$|\rho(j)-\pi(j)|\le n-1$, and $|k_jm|\ge m> n$, if $k_j\ne 0$, so   
if at least one $k_j \in \Z$ is different from 0, then there exists no solution. This proves
the theorem.
\end{biz}\bv

\section{Asymptotics of the trace}\label{masodik}

Let  $U(n)=(U(n)_{ij}):\Omega\rightarrow \mathcal U(n)$ be a Haar distributed 
unitary random matrix. In this section we are interested in the
convergence of $\tr U(n)$ as $n \to \infty$. Since the correlation 
between the diagonal entries decreases with $n$, one expects on the 
basis of the central limit theorem, that the limit of the trace has 
complex normal distribution. We prove this by the method of moments.

\begin{tetel}\label{1hatv}
Let $U(n)$ be a sequence of $n\times n$ Haar unitary random matrices.
Then $\tr U(n)$ converges in distribution to a standard complex normal
random variable as $n \to \infty$.
\end{tetel} 

\begin{biz}
For the sake of simplicity we write $U$ instead of $U(n)$.
First we study the asymptotics of the moments
\begin{eqnarray*}
\lefteqn{E\left((\tr U)^k(\overline{\tr U})^k\right)=
E\Big(\Big(\sum_{i=1}^n U_{ii}\Big)^k
\Big(\sum_{j=1}^n \overline U_{jj}\Big)^k\Big)}\\
&&=\sum_{i_1,\dots ,i_k=1}^n\sum_{j_1,\dots,j_k=1}^n E(U_{i_1i_1}\dots U_{i_ki_k}
\overline U_{j_1j_1}\dots \overline U_{j_kj_k}),\label{sum}
\end{eqnarray*}
$k\in \Z^+$. By Lemma \ref{mom0} parts of the above sum are zero, we
need to consider only those sets of indices $\{i_1,\dots,i_k\}$ and 
$\{j_1,\dots,j_k\}$ which coincide (with multiplicities). Look at a summand 
$E(|U_{i_1i_1}|^{2k_1}\dots |U_{i_ri_r}|^{2k_r})$, where $\sum_{l=1}^r k_l=k$. 
>From the  H\"older inequality 
\begin{equation}\label{kismom}
E(|U_{i_1j_1}|^{2k_1}\dots |U_{i_rj_r}|^{2k_r})\le\prod_{l=1}^r
\sqrt[2^l]{E(|U_{i_lj_l}|^{2\cdot 2^lk_l})}=\prod_{l=1}^r
{{n+2^lk_l-1}\choose {2^lk_l-1}}^{-1/2^l}= O\left(n^{-k}\right).
\end{equation}
The number of those sets of indices, where among the numbers $i_1,\dots,i_k$
there are at least two equal is at most
$$
k!{k\choose 2}n^{k-1}= O(n^{k-1}).
$$
By (\ref{kismom}) the order of magnitude of these factors is
$O(n^{-k})$, so this part of the sum tends to zero as $n \rightarrow \infty$.

Next we assume that $i_1,\dots,i_k$ are different. Since 
by translation invariance any row or column can be replaced by any other, we have
\begin{equation}\label{mk}
E(|U_{i_1i_1}|^{2}\dots |U_{i_ki_k}|^{2})=
E(|U_{11}|^{2}\dots |U_{kk}|^{2}))=:M_k^n.
\end{equation}
It is enough to determine this quantity and to count how many of these
terms are in the trace.

The length of the row vectors of the unitary matrix is 1, hence 
\begin{equation}\label{szum}
\sum_{i_1=1}^n \dots \sum_{i_k=1}^n E\left(|U_{i_11}|^2\dots |U_{i_kk}|^2\right)=1.
\end{equation}
We divide the sum into two parts: the number of terms with different
indices is ${n!/(n-k)!}$, and again the
translation invariance implies that each of them equals to $M_k^n$, and we denote by
$\varepsilon_k^n$ the sum of the other terms. Therefore
\[
\varepsilon_k^n=1-\frac{n!}{(n-k)!}M_k^n\le k!{k\choose 2} O(n^{-k})\to 0,
\]
and
\[M_k^n=\frac{(1-\varepsilon_k^n)(n-k)!}{n!}.\]
Now we can count how many expectations of value $M_k^n$ are there in the sum 
(\ref{sum}). We can fix 
the indices $i_1,\dots,i_k$ in ${n!}/{(n-k)!}$ ways, and
we can permute them in $k!$ ways to get the indices 
$j_1,\dots, j_k$. The obtained equation
\[\lim_{n\rightarrow\infty} E\left((\tr U_n)^k(\overline{\tr U_n})^k\right)=
\lim_{n\rightarrow\infty}\frac{n!}{(n-k)!}k!\frac{(1-\varepsilon_k^n)(n-k)!}{n!}
=k!
\]
finishes the proof.

For the mixed moments we have by Lemma \ref{mom0}
$$
E\left((\tr U_n)^k(\overline{\tr U_n})^m\right)=0 
\qquad (k\ne m),
$$
and we have proven the convergence of all moments. The only thing is
left to conclude the convergence in distribution is to show that the
moments determine uniquely the limiting distribution ( VIII. 6 in 
\cite{feller2}).
By the Carleman criterion (in VII. 3 of \cite{feller2}) for a real valued random variable the  
moments $\gamma_k$ determine the distribution uniquely if 
\[
\sum_{k \in \N} \gamma_{2k}^{-\frac 1k}=\infty.
\]
Although we have complex random variables, the distribution of the argument
is uniform, and we can consider them as real valued random variables. 
The Stirling formula tells us
\[
\sum_{k \in \N} (k!)^{-\frac 1k}\ge \sum_{k\ge M}
\left(\left(\frac {2k}{\textrm{e}}\right)^k
\right)^{-\frac 1k}= \frac{\textrm{e}}{2}\sum_{k\ge M}\frac 1k=\infty.
\] 
for a large $M\in \N$, since $\sqrt{2k\pi}\le 2^k$, if $k\ge 2$.
\end{biz}\bv

\section{Asymptotic behaviour of the traces of higher powers}\label{harmadik}

The aim of this section is to study the trace of higher powers of a Haar unitary. 
This was done also by Diaconis and Shashahani in \cite{diaconis}. Here we use
elementary methods.

\begin{tetel}
Let $Z$ be standard complex normal distributed \vv, 
then for the sequence of $U_n$ $n\times n$ Haar unitary random matrices
$\tr U_n^l \rightarrow \sqrt l Z$ in distribution.
\end{tetel}

\begin{biz}
We use the method of moments again. Lemma \ref{mom0}  implies that 
we only have to take into consideration
$E\left(\left(\tr U_n^l\right)^k\left(\overline{\tr U_n^l}\right)^k\right)$, for all 
$k\in\Z^+$. 
\begin{eqnarray*}
\lefteqn{E\left(\left(\tr U_n^l\right)^k\left(\overline{\tr U_n^l}\right)^k\right)}\\
&&\kern-.25in=E\Bigg(\Big(\sum_{i_1,\dots,i_l}
U_{i_1i_2}U_{i_2i_3}\dots U_{i_{l-1}i_l}U_{i_li_1}
\Big)^k\Big(\sum_{j_1,\dots, j_l}\overline U_{j_1j_2}\overline U_{j_2j_3}
\dots \overline U_{j_{l-1}j_l}\overline U_{j_lj_1}\Big)^k\Bigg)\\ \\
&&\kern-.25in=\sum E\left(U_{i_1i_2}\dots U_{i_li_1}U_{i_{l+1}i_{l+2}}\dots
U_{i_{2l}i_{l+1}}\dots U_{i_{l(k-1)+1}i_{l(k-1)+2}}\dots
U_{i_{kl}i_{l(k-1)+1}}\right.\\
&&\kern-.05in\times \left.\overline U_{j_1j_2}\dots \overline U_{j_lj_1}\overline 
U_{j_{l+1}j_{l+2}}\dots
\overline U_{j_{2l}j_{l+1}}\dots\overline U_{j_{l(k-1)+1}j_{l(k-1)+2}}\dots
\overline U_{j_{kl}j_{l(k-1)+1}}\right),
\end{eqnarray*}
where the indices $i_1,\dots, i_{kl},j_1,\dots, j_{kl}$ run from 1 to $n$, and
by Lemma \ref{mom0} if the sets $\{i_1,\dots , i_{kl}\}$ and $\{j_1,\dots, j_{kl}\}$
are different, then the expectation of the product is zero.
It follows from the the Cauchy and H\"older inequalities, and (\ref{kismom}), that
\begin{eqnarray}\label{pici}
\lefteqn{\left|E\left(U_{i_1i_2}\dots U_{i_{kl}i_{l(k-1)+1}}\overline U_{j_1j_2}
\dots \overline U_{j_{kl}j_{l(k-1)+1}}\right)\right|}\nonumber\\ 
&&\le E\left|U_{i_1i_2}\dots U_{i_{kl}i_{l(k-1)+1}}\overline U_{j_1j_2}
\dots \overline U_{j_{kl}j_{l(k-1)+1}}\right|\\
&& \le\sqrt{E\left(|U_{i_1i_2}|^2\dots |U_{i_{kl}i_{l(k-1)+1}}|^2|\overline U_{j_1j_2}|^2
\dots |\overline U_{j_{kl}j_{l(k-1)+1}}|^2\right)}\le O\left(n^{-kl}\right). \nonumber
\end{eqnarray}
Again the number of the set of indices, where there exist at least two equal
indices is at most $O(n^{kl-1})$, so the sum of the corresponding
expectations tends to zero as $n\rightarrow\infty$.
Suppose that all the indices are different. There exist 
$\frac{n!}{(n-kl)!}(kl)!=O(n^{kl})$
of these kinds of index sets, and now we will prove, that most of the 
corresponding products
have order of magnitude less than $n^{-kl-1}$. Consider for any $0\le r\le kl$
\[N_k^n(r):=E\left(|U_{12}|^2|U_{23}|^2\dots |U_{r1}|^2U_{r+1,r+2}\dots
U_{kl-1,kl}U_{kl,r+1}\overline U_{r+2,r+1}\dots \overline U_{r+1,kl}\right) .\]
Note that $N_k^n(kl)=N_k^n(kl-1)=M_{kl}^n$, and if $\{i_1,\dots
i_{kl}\}=\{j_1,\dots ,j_{kl}\}$, and all the 
indices are different, then the corresponding term equals to $N_k^n(r)$ for 
some $0\le r\le kl$.
Using the orthogonality of the rows for $0\le r\le kl-2$
\begin{equation}\label{orto}
E\left(\sum_{j=1}^n|U_{12}|^2|U_{23}|^2\dots |U_{r1}|^2U_{r+1,r+2}\dots
U_{kl-1,j}U_{kl,r+1}\overline U_{r+2,r+1}\dots \overline U_{r+1,j}\right)=0.
\end{equation}
If $j\ge kl$, then the permutation invariance implies, that
\[
E\left(|U_{12}|^2|U_{23}|^2\dots |U_{r1}|^2U_{r+1,r+2}\dots
U_{kl-1,j}U_{kl,r+1}\overline U_{r+2,r+1}\dots \overline U_{r+1,j}\right)=N_k^n(r),\]
so we can write from (\ref{orto}) 
\[
(n-kl)N_k^n(r)=-E\left(\sum_{j=1}^{kl}|U_{12}|^2|U_{23}|^2\dots 
|U_{r1}|^2U_{r+1,r+2}\dots
U_{kl-1,j}U_{kl,r+1}\overline U_{r+2,r+1}\dots \overline U_{r+1,j}\right).
\]
On the right side there is a sum of $kl$ numbers which are less than $O(n^{-kl})$
because of (\ref{pici}), so this equation holds only if $N_k^n(r)\le O(n^{-kl-1})$. 

\bigskip\noindent 
We have to compute the sum of the expectations
\[E\left(|U_{i_1i_2}|^2\dots|U_{i_li_1}|^2\dots|U_{i_{(k-1)l+1}i_{(k-1)l+2}}|^2
\dots |U_{i_{kl}i_{(k-1)l+1}}|^2\right)=M_{kl}^n.\]
Now we count the number of these summands, so first we fix the set of sequences 
of length $l$
\sloppy${I_{l,k}=\{(i_{(u-1)l+1},\dots,i_{ul}),1\le u \le k\}}$, and we try to find 
the set \sloppy${J_{l,k}=\{(j_{(u-1)l+1},\dots,j_{ul}),1\le u \le k\}}$, 
which gives $M_{kl}^n$. If the product
contains $U_{i_ri_{r+1}}$, then it has to contain $\overline U_{i_ri_{r+1}}$, so if
$i_r$ and $i_{r+1}$ are in the same sequence of $I_{l,k}$, then $j_s=i_r$ and $j_t=i_{r+1}$
have to be in the same sequence of $J_{l,k}$, and $t=s+1$ modulo $l$. This means, that
for all $1\le u\le k$  there exists 
a sequence $(i_{(v-1)l+1},\dots, i_{vl})\in I_{k,l}$ and a cyclic permutation $\pi$ of 
the numbers $\{(v-1)l+1,\dots, vl\}$ such that
$(j_{(u-1)l+1},\dots, j_{ul})=(i_{\pi((v-1)l+1)},\dots,i_{\pi(vl)})$. 
We conclude, that for each $I_{l,k}$ there are $k!l^k$ $J_{l,k}$, since we can permute
the sets of $I_{l,k}$ in $k!$ ways, and in all sets there are $l$ cyclic permutations.
Clearly there are $\frac{n!}{(n-kl)!}$ sets $I_{l,k}$, so 
\[\lim_{n\rightarrow\infty}E\left(\left(\tr U_n^l\right)^k\left(\overline{\tr U_n^l}\right)^k
\right)=\lim_{n\rightarrow\infty}\frac{n!}{(n-kl)!}k!l^k
\frac{(1-\varepsilon_{kl}^n)(n-kl)!}{n!}=k!l^k,\]
and as in the proof of Theorem \ref{1hatv} this is the $k$th moment of $(\sqrt l Z)\overline{(\sqrt lZ)}$.
\end{biz}\bv

\section{Independence of the trace of the powers}\label{negyedik}

In this section we prove that the limits of the trace of different powers 
are independent. The method of computation is the same as in the previous sections.

\begin{tetel}
Let $U_n$ be a sequence of Haar unitary random matrices as above. Then $\tr U_n,\tr 
U_n^2,\dots,
\tr U_n^l$ are asymptotically independent.
\end{tetel}

\begin{biz}
We will show, that the joint moments of $\tr U_n, \tr U_n^2, \dots \tr U_n^l$ converge to the
joint moments of $Z_1,  \sqrt 2 Z_2,\dots, \sqrt l Z_l$, where $Z_1, Z_2, \dots Z_l$
are independent standard complex normal \vv s. The latter joint moments are 
\[E\left(\prod_{i=1}^li^{\frac{a_i+b_i}{2}}Z_i^{a_i}\overline Z_i^{b_i}\right)=
\prod_{i=1}^li^{\frac{a_i+b_i}{2}}E\left(Z_i^{a_i}\overline Z_i^{b_i}\right)=
\prod_{i=1}^l\delta_{a_ib_i}a_i!i^{a_i}.\]
>From Lemma \ref{mom0}, if $\sum_{i=1}^lia_i\ne \sum_{i=1}^lib_i$, then the moment
\[E \left(\prod_{i=1}^l\left(\tr U_n^i\right)^{a_i}\left(\overline{\tr U_n^i}\right)^{b_i}\right)=0.\]
Now $\sum ia_i= \sum ib_i$. Again if the indices in the first and the second part are 
not the same, then the expectation is zero according to Lemma $\ref{mom0}$.
The order of magnitude of each summand is at most $O\left(n^{-\sum ia_i}\right)$, as above,
so if not all the indices are different, then the sum of these expectations tends to zero, as
$n \rightarrow \infty$. The same way as in the proof of the previous theorem, those summands
where there is a $U_{i_ri_{r+1}}\overline U_{i_ri_s}$, $i_{r+1}\ne i_s$ in the product
are small.  So now we have to sum the expectations $M_{\sum ia_i}^n$. If we fix
the set of first indices $I$, then again the sequences of the appropriate
$J$, have to be cyclic permutations of the sequences of $I$. So the number of 
the sequences of length $i$ in $I$ is the same as in $J$, which means
$a_i=b_i$ for all $1\le i\le l$. 
The number of the $I$ sets is
$\frac{n!}{\left(n-\sum ia_i\right)!}$, so we have arrived to
\begin{eqnarray*}
\lefteqn{\lim_{n\rightarrow\infty} E \left(\prod_{i=1}^l\left(\tr U_n^i\right)^{a_i}
\left(\overline{\tr U_n^i}\right)^{b_i}\right)}\\
&&=\lim_{n\rightarrow \infty}\frac{n!}{(n-\sum ia_i)!}\prod_{i=1}^l\delta_{a_i,b_i}i^{a_i}a_i!
\frac{\left(1-\varepsilon_{\sum ia_i}^n\right)\left(n-\sum ia_i\right)!}{n!}=
\prod_{i=1}^l\delta_{a_i,b_i}a_i!i^{a_i}.
\end{eqnarray*}
\end{biz}\bv

\section{Truncation}

Let $U$ be an $n \times n$ Haar distributed unitary matrix. By truncating $n-m$
bottom rows and $n-m$ last columns, we get an $m \times m$ matrix $U_{[n,m]}$.
In this section we study the limit of $U_{[n,m]}$ when $n \to \infty $
and $m$ is fixed. Our method is based on the explicite form of the joint
eigenvalue density.

The truncated matrix is not unitary but it is a contraction. Hence the
eigenvalues $z_1 , z_2 , \ldots , z_m \in D^m$, where $D = \{ z \in \C \, :
\, |z| \leq 1\}$ is the unit disc. According to \cite{ZS} the joint
probability density of the eigenvalues is
$$
C_{[n,m]} \prod_{i < j} |\zeta_i -\zeta_j |^2 \prod_{i=1}^m (1-|\zeta_i |^2
)^{n-m-1}
$$
on $D^m$. 

Since the normalizing constant $C_{[n,m]}$ was not given in \cite{ZS},
we first compute it by integration. To do this, we write $\zeta_i =r_i
e^{\im \varphi_i}$ and $d \zeta_i = r_i \, d r_i \, d \varphi_i $. Then
\begin{eqnarray*}
\lefteqn{C_{[n,m]}^{-1}
=\int_{D^m}\prod_{1\le i<j\le m}|z_i-z_j|^2\prod_{i=1}^m(1-|z_i|^2)^{n-m-1}\, dz}\\
&&=\int_{[0,1]^m}\int_{[0,2\pi]^m}\prod_{1\le i<j\le m}|r_ie^{i\varphi_1}-r_je^{i\varphi_j}|^2
\prod_{i=1}^m(1-r_i^2)^{n-m-1}\prod_{i=1}^m r_i\ud \varphi\ud r.
\end{eqnarray*}
Next we integrate with respect to 
$\ud\varphi=\ud \varphi_1\ud \varphi_2\dots \ud \varphi_m$  by
transformation into complex contour integral what we evaluate by means of the
residue theorem.
\begin{eqnarray*}
\int_{[0,2\pi]^n}\prod_{1\le i<j\le m}|r_ie^{i\varphi_1}-r_je^{i\varphi_j}|^2
\ud \varphi&&\kern-.24in \\
&&\kern-1.9in=(-i)^m\int_{\T^n}\prod_{1\le i<j\le m}|r_iz_i-r_jz_j|^2
\prod_{i=1}^mz_i^{-1}\ud z \\
&&\kern-1.9in=(-i)^m\int_{\T^n}\prod_{1\le i<j\le m}(r_iz_i-r_jz_j)
(r_iz_i^{-1}-r_jz_j^{-1})
\prod_{i=1}^mz_i^{-1}\ud z \\
&&\kern-1.9in=(-i)^m\int_{\T^n}\prod_{i=1}^mz_i^{-1}\det
\left[\begin{array}{llll}
1&1&\dots&1\\
r_1z_1^{}&r_2z_2^{}&\dots&r_mz_m^{}\\
\vdots&&\ddots&\vdots\\
r_1^{m-1}z_1^{m-1}&r_2^{m-1}z_2^{m-1}&\dots&r_m^{m-1}z_m^{m-1}
\end{array}\right]\times\\
&&\kern-1.11in\times
\det\left[\begin{array}{llll}
1&1&\dots&1\\
r_1z_1^{-1}&r_2z_2^{-1}&\dots&r_mz_m^{-1}\\
\vdots&&\ddots&\vdots\\
r_1^{m-1}z_1^{-(m-1)}&r_2^{m-1}z_2^{-(m-1)}&\dots&r_m^{m-1}z_m^{-(m-1)}
\end{array}\right]
\ud z \\
&&\kern-1.9in=(-i)^m\int_{\T^n}\prod_{i=1}^mz_i^{-1}\sum_{\pi \in S_m}(-1)^{\sigma
(\pi)}\prod_{i=1}^m(r_iz_i)^{\pi(i)-1}
\sum_{\rho \in S_m}(-1)^{\sigma(\rho)}\prod_{i=1}^m(r_iz_i^{-1})^{\rho(i)-1}
\ud z \,.
\end{eqnarray*}
We have to find the coefficient of $\prod_{i=1}^mz_i^{-1}$, this gives
that only $\rho=\pi$ contribute and the integral is
$$
(2\pi)^m\sum_{\rho \in S_m}\prod_{i=1}^m(r_i)^{2(\rho(i)-1)}.
$$
So we have 
\begin{eqnarray*}
C_{[n,m]}^{-1}&&\kern-.24in=(2\pi)^m\int_{[0,1]^m}\sum_{\rho \in S_m}
\prod_{i=1}^m(r_i)^{2(\rho(i)-1)}
\prod_{i=1}^m(1-r_i^2)^{n-m-1}\prod_{i=1}^m r_i \ud r\\
&&\kern-.24in =(2\pi)^mm!\prod_{i=1}^m\int_0^1 r_i^{2i-1}(1-r_i^2)^{n-m-1} \ud r_i
\end{eqnarray*}
and the rest is done by integration by parts:
\begin{eqnarray*}
\lefteqn{\int_0^1 r^{2k+1}(1-r^2)^{n-m-1}\ud r=\frac{k}{n-m}\int_0^1
r^{2k-1}(1-r^2)^{n-m}\ud r}\\
&&=\frac{k!}{(n-m)\dots (n-m+k-1)}\int_0^1 r(1-r^2)^{n-m+k-1}\ud r\\
&&={n-m+k-1\choose k}^{-1}\frac{1}{2(n-m+k)}.
\end{eqnarray*}
Therefore 
\[
C_{[n,m]}^{-1}= \pi^mm! \prod_{k=0}^{m-1}{n-m+k-1\choose k}^{-1}
\frac{1}{n-m+k}.
\]

Now we consider  $\sqrt{n/m} \,U_{[n,m]}$. Its joint probability dnsity
of the eigenvalues is simply derived from the above density of $U_{[n,m]}$ by the transformation
\[(\zeta_1,\dots ,\zeta_m)\mapsto \left(\sqrt\frac mn\zeta_1,\dots ,\sqrt \frac mn \zeta_m
\right),\]
and it is given as
\begin{eqnarray*}
\lefteqn{C_{[n,m]}\left(\frac mn\right)^m\prod_{i<j} \left|\sqrt{\frac mn}\zeta_i-\sqrt{\frac mn}
\zeta_j\right|^2\prod_{i=1}^m \left(1-\frac{m|\zeta_i|^2}{n}\right)^{n-m-1}}\\
&&=\frac{1}{\pi^mm!}\prod_{k=0}^{m-1}{{n-m+k-1}\choose k}(n-m+k)\left(
\frac mn\right)^{m(m+1)/2}\\
&&\kern1in\times\prod_{i<j}|\zeta_i-\zeta_j|^2\prod_{i=1}^m\left(1-\frac{m|\zeta_i|^2}{n}\right)^{n-m-1}\\
&&=\frac{1}{\pi^mm!}\prod_{k=0}^{m-1}\frac{n^{k+1}(1+o(1))}{k!}\left(
\frac mn\right)^{m(m+1)/2}\prod_{i<j}|\zeta_i-\zeta_j|^2\prod_{i=1}^m\left(1-\frac{m|\zeta_i|^2}{n}
\right)^{n-m-1}\\
&&=\frac{m^{m(m+1)/2}}{\pi^m\prod_{k=1}^m k!}(1+o(1))\prod_{i<j}|\zeta_i-\zeta_j|^2\prod_{i=1}^m
\left(1-\frac{m|\zeta_i|^2}{n}\right)^{n-m-1}.
\end{eqnarray*}
The limit of the above as $n\to \infty$ is
\[\frac{m^{m(m+1)/2}}{\pi^m\prod_{k=1}^m k!}\exp\left(-m\sum_{i=1}^m|\zeta_i|^2\right)\prod_{i<j}|\zeta_i-\zeta_j|^2,\]
which is exactly the joint eigenvalue density of the standard $m\times m$ non-selfadjoint Gaussian matrix.
This implies the following theorem.

\begin{tetel}
The normalized truncated matrix 
$$
\sqrt{\frac nm} U_{[n,m]}
$$
converge in distribution to the standard $m\times m$ non-selfadjoint
Gaussian matrix.
\end{tetel}

Details of this convergence will be subject of a forthcoming publication.

\end{document}